\documentclass{ifacmtg}
\usepackage{amsfonts,amssymb,latexsym,amsmath,amscd,theorem,epic,epsfig,graphics}
\newtheorem{theorem}{Theorem}
\newtheorem{remark}{Remark}
\newtheorem{lemma}{Lemma}

\begin{document}
\newcommand{\eps}{{\varepsilon}}
\newcommand{\proofend}{$\Box$\bigskip}
\newcommand{\sign}{{\mbox{\rm sign}}}
\newcommand{\mes}{{\mbox{\rm mes}}}
\newcommand{\lcm}{{\mbox{\rm lcm}}}
\newcommand{\Id}{{\mbox{\rm Id}}}
\newcommand{\tet}{{\theta}}
\newcommand{\CC}{{\mathbb C}}
\newcommand{\QQ}{{\mathbb Q}}
\newcommand{\R}{{\mathbb R}}
\newcommand{\ZZ}{{\mathbb Z}}
\newcommand{\PP}{{\mathbb P}}
\newcommand{\NN}{{\mathbb N}}
\newcommand{\Del}{{\Delta}}
\newcommand{\dspace}{\lineskip=2pt\baselineskip=18pt\lineskiplimit=0pt}
\newcommand{\sspace}{\lineskip=2pt\baselineskip=14pt\lineskiplimit=0pt}
\newcommand{\bet}{{\beta}}
\newcommand{\oa}{{\overline a}}
\newcommand{\un}{{\bf n}}
\newcommand{\ob}{{\overline b}}
\newcommand{\oc}{{\overline c}}
\newcommand{\um}{{\bf m}}
\newcommand{\uq}{{\bf q}}
\newcommand{\ui}{{\bf i}}
\newcommand{\uj}{{\bf j}}
\newcommand{\kap}{{\kappa}}
\newcommand{\del}{{\delta}}
\newcommand{\sig}{{\sigma}}
\newcommand{\alp}{{\alpha}}
\newcommand{\Sig}{{\Sigma}}
\newcommand{\Gam}{{\Gamma}}
\newcommand{\gam}{{\gamma}}
\newcommand{\Lam}{{\Lambda}}
\newcommand{\lam}{{\lambda}}
\newcommand{\const}{{\operatorname{const}}}
\runauthor{Shustin and Fridman}
\begin{frontmatter}
\title
{ An Improved Stability Method for Linear Systems with
Fast-Varying Delays}
\author[Paestum]{Eugenii Shustin}
\author[Paestum1]{Emilia Fridman}
\address[Paestum]{
School of Mathematics, Tel Aviv University, Tel Aviv 69978,
Israel\\
e-mail: shustin@post.tau.ac.il }
\address[Paestum1]{
School of Electrical Eng. - Systems, Tel Aviv University, Tel Aviv
69978,
Israel\\
e-mail: emilia@eng.tau.ac.il }

\begin{abstract}
Stability of linear systems  with uncertain bounded time-varying
delays is studied under assumption that the nominal delay values
are not equal to zero.
 An
input-output approach to stability of such systems  is known to be
based on the bound of the $L_2$-norm of a certain integral
operator. There exists a bound on this operator in two cases: in
the case where the delay derivative is not greater than $1$ and in
the case without any constraints on the delay derivative. In the
present note we fill the  gap between the two cases by deriving a
tight operator  bound which  is an increasing and continuous
function of the delay derivative  upper bound  $d \geq 1$. For
$d\to \infty$ the new bound  corresponds to the second case and
improves the existing bound. As a result,
delay-derivative-dependent frequency-domain and time-domain
stability criteria are derived for systems with the delay
derivative greater than $1$.
\end{abstract}
\begin{keyword}
 time-varying delay, stability, input-output
approach, $L_2$-norm.
\end{keyword}
\end{frontmatter}

\section{Introduction}

 Two main approaches have been applied to stability analysis of
linear systems with uncertain time-varying delay: a direct
Lyapunov approach and an input-output approach
 (see e.g. Gu {\it et al.}, 2003),
 which reduces the stability
analysis of the uncertain system to the analysis of the class of
systems with the same nominal part but with additional inputs and
outputs.
In the existing literature the uncertain
 time-varying delay
 has been divided into
 two types: the {\it slowly-varying} delay
 (with  delay derivative less than $d<1$) and the
 {\it fast-varying delay} (without any constraints on the delay
 derivative) (see e.g.
Kolmanovskii \& Myshkis , 1999; Niculescu , 2001).
 Recently a third type of  {\it moderately
varying delay} has been revealed in Fridman \& Shaked (2005),
where the delay derivative is not greater than $1$ (almost for all
$t$). This has been obtained by applying the input-output approach
to stability. It is known Gu {\it et al.} (2003), Kao \&  Lincoln
(2004) that the latter approach to systems with time-varying
bounded delays is based on the bound of the $L_2$-norm of a
certain integral operator.

In the present paper we fill the  gap between the case of the
delay derivative not greater than $1$ and the fast-varying delay
by deriving a new integral operator bound. This bound is an
increasing and continuous function of the delay derivative bound
$d \geq 1$. In the limit case (where $d\to \infty$) which
corresponds to the fast-varying delay, the new bound improves the
existing one. As a result, improved
frequency-domain and time-domain stability criteria are derived
for systems with the delay derivative greater than $1$.

{\bf Notation:}\ Throughout the paper the superscript `$T$' stands
for matrix transposition, ${\cal R}^n$ denotes the $n$ dimensional
Euclidean space with vector norm $\|\cdot \|$, ${\cal R}^{n\times
m}$ is the set of all $n\times m$ real matrices, and the notation
$P\!>\!0$,
 for  $P\in{\cal R}^{n\times n}$
means that  $P$ is symmetric and positive  definite.
 The
symmetric elements of the symmetric matrix will be denoted by
${*}$.
$L_2$ is the space of square integrable functions $v:[0,\infty)\to
C^n$  with the norm $ \|v\|_{L_2} = [\int_0^{\infty}
\|v(t)\|^2dt]^{1/2},$  $\|A\|$ denotes the Euclidean norm of a
$n\times n$ (real or complex) matrix $A$, which  is equal to the
maximum singular value of $A$.
 For a
  transfer function  matrix of a stable system $
 G(s)$, $s\in C$
$$\|G\|_{\infty}=\sup_{-\infty <w<\infty}\|G(iw)\|, \quad i=\sqrt{-1}.$$

\section{Problem Formulation}
We consider the following linear system with  uncertain
time-varying delay $\tau(t)$ :
\begin{equation}
\label{x}
   \dot x(t)=A_0x(t)+A_1 x(t-\tau (t)),
\end{equation}
 where $x(t)\in{\cal R}^{n} $ is the system
state,  $A_i, \ i=0,1$  are constant matrices.

The uncertain delay $\tau(t)$ has a form
\begin{equation}
\label{tau} \tau(t)=h+\eta(t),  \  \quad
 |\eta(t)|\leq \mu \leq
h,
\end{equation}
  where $h$ is a   known nominal delay value   and $\mu$ is a  known
  upper bound on the delay uncertainty. In the existing literature
  Kolmanovskii \& Myshkis (1999), Niculescu (2001), Gu {\it et al.} (2003) 
  the following types of uncertain time-varying delays are usually
  considered:

{\bf Case A (slowly-varying delay):}  $\tau(t)$ is a
differentiable almost everywhere function, satisfying
\begin{equation}
\label{h} \dot \tau (t)=\dot \eta (t) \leq d=1+p,
\end{equation}
where $-1\leq p<0$;

{\bf Case B (fast-varying delay):}  $\tau(t)$ is a measurable
(e.g. piecewise-continuous) function.$ \vspace{0.2cm}$

Recently a {\it moderately-varying} delay with $\dot \tau (t)\le d
= 1$ was introduced in Fridman \& Shaked (2005). In the present
note we enlarge the latter class of delays as follows:

{\bf Case C (moderately-varying delay):}   $\tau(t)$ is a
differentiable almost everywhere function, satisfying (\ref{h})
with $ p\geq 0$.

\vspace{.3cm}

In the present note we will improve the stability results in cases
B and C by applying input-output approach and by deriving new
inequalities. The results are easily generalized to the case of
any finite number of the delays.

We represent (\ref{x}) in the form:
 \begin{equation}
 \begin{array}{lll}
 \label{x1}
 \dot x(t)=A_0x(t)+A_1x(t-h)- A_1\int_{
 -h-\eta}^{-h}\dot x(t+s)ds.
 \end{array}
 \end{equation}
Following Fridman \& Shaked (2005) we introduce the following
auxiliary system:
 \begin{equation}
 \label{input}
\begin{array}{lll}
\dot x(t)=A_0x(t)+ A_1x(t-h)+
 \mu A_1
 u(t),
 \\
 y(t)=\sqrt{ {\cal F}(p)}\dot x(t), \\
\end{array}
\end{equation}
 with the feedback
\begin{equation}
\label{output}
\begin{array}{lll}
u(t)=-\frac{1}{\mu\cdot\sqrt {{\cal F}(p)}}\int_{-h-
 \eta}^{-h } y(t+s)ds,
\end{array}
\end{equation}
where ${\cal F}: [-1, \infty] \to R^{+}$ is a scalar function
which will be shortly defined and $p$ is given by (\ref{h}). The
results for the delay of case B correspond to $p=\infty$, i.e. to
${\cal F (\infty})$ in the input-output model (\ref{input}),
(\ref{output}). Substitution of (\ref{output}) in (\ref{input})
readily leads to (\ref{x1}).

We are looking for ${\cal F}(p)$ which satisfies the following
inequality
\begin{equation}
\label{u} ||u||_{L_2}^2\le ||y||_{L_2}^2, \quad \forall\ y\in
L_2[0 ,\infty),\ y\big|_{[-\infty,0]}\equiv 0,
\end{equation}
where $u$ is given by (\ref{output}). This is equivalent to the
fact that $\mu {\sqrt {{\cal F}(p)}}$ is an upper bound on the
$L_2$-norm of the integral operator $\Delta:L_2[0,\infty)\to
L_2[0,\infty)$
\begin{equation}
\label{Del} z(t)=\Delta y (t) = \int_{-h-
 \eta}^{-h } y(t+s)ds,\ y\big|_{[-\infty,0]}\equiv 0,
 \end{equation}
 i.e. that
 \begin{equation}\begin{array}{lll}
 \|z\|_{L_2}^2\le\mu^2{\cal F}(p)\|y\|_{L_2}^2\ ,\\
  \forall y\in L_2[0
,\infty),\ y\big|_{[-\infty,0]}\equiv 0. \label{en1}
 \end{array}
 \end{equation}
Our objective is to find ${\cal F}(p)$ (as small as possible) such
that (\ref{u}) (or equivalently (\ref{en1})) holds.

For $-1\le p<0$ (case A) it was established in Gu {\it et al.}
(2003) that ${\cal F}(p)$ can be chosen to be $1$.
 For $p\ge 0$ the following was found in
Fridman \& Shaked (2005):  ${\cal F}(0)=1$ and ${\cal F}(p)\equiv
2$ for $p\in (0, \infty]$.

We note that the value $1$ of ${\cal F}(p)$ for $-1\le p\le 0$ can
not be improved (i.e. chosen to be less than $1$). Indeed, taking
constant delay $\eta(t)\equiv\mu$, which satisfies the condition
of case A for any $-1\le p\le 0$,
 we consider the functions
$y_\theta(t)=1$ as $0\le t\le\theta$, and $y_\theta(t)=0$ as
$t>\theta$. Using formula (\ref{output}) with ${\cal F}(p)=1$ we
immediately obtain
$$\|y_\theta\|_{L_2}^2=\theta^2,\quad \|u\|_{L_2}^2=(\theta-\mu)^2+\frac{2}{3}\mu^2\ ,$$
and hence
$\|u\|_{L_2}/\|y_\theta\|_{L_2}\to 1$ as $\theta\to\infty$.

In the present paper we will improve the values of ${\cal F}(p)$
for $p>0$ by showing that ${\cal F}(p)$ can be chosen as a {\it
continuous increasing function of $p\geq 0$} satisfying ${\cal
F}(0)=1$ (as in Fridman \& Shaked (2005)), but ${\cal F}(p)< {\cal
F}(\infty)= 1.75$ for $p>0$. The improved values of  ${\cal F}(p)$
will readily lead to improved stability criteria.

\section{Main Results}

\subsection{New Bounds}

Proofs of the Lemmas of this section are given in the Appendix.

\begin{lemma}\label{l2} Consider case C. For all $y(t)\in L_2[0, \infty)$ and such that
$y(t)=0  \ \forall t\leq 0$ and for $u(t)$ given by (\ref{output})
inequality (\ref{u}) holds with ${\cal F}$ given by
\begin{equation}
\label{lem2}
{\cal F}(p)=\begin{cases}\frac{2p+1}{p+1},\quad&\text{if}\quad 0\leq p < 1,\\
\frac{7p-1}{4p},\quad&\text{if}\quad p\geq 1 \ .\end{cases}
\end{equation}
\end{lemma}

As it was mentioned above, ${\cal F} $ is increasing continuous
function satisfying  for $p>0$ the following inequality: $1={\cal
F}(0)<{\cal F}(p)< \lim_{p\to \infty} {\cal F}(p) = 7/4.$

\begin{lemma}\label{l1} Consider case B. For all $y(t)\in L_2[0, \infty)$ and such that
$y(t)=0  \ \forall t\leq 0$ and for $u(t)$ given by (\ref{output})
 inequality (\ref{u}) holds with ${\cal F}(\infty):={7/ 4}$.
\end{lemma}

\begin{remark} The value $7/4=1.75$ for ${\cal F}(\infty)$ in Lemma \ref{l1}
is not far from an optimal one. The following example shows that
it cannot be less than $1.5$. Namely, define scalar functions
$y(t)$ and $\eta(t)$ by
$$y(t)=\begin{cases}t,\quad&\text{if}\ 0\le t\le\mu,\\
\mu-t,\quad&\text{if}\ \mu\le t\le 2\mu,\\ 0,\quad&\text{if}\
y(2\mu-y)<0,\end{cases}$$
$$\eta(t)=\begin{cases}-\mu,\quad&\text{if}\ t\le\mu,\\
\mu,\quad&\text{if}\ t>\mu.\end{cases}$$ Setting in (\ref{output})
${\cal F}(\infty )=3/2$ we have
$u(t)=-\frac{1}{\mu\sqrt{3/2}}z(t)$, where
$$z(t+h)=\int_{t-\eta(t)}^ty(s)ds$$
$$=
\begin{cases}-(t+\mu)^2/2,\quad&\text{if}\
-\mu\le t+h\le 0,\\ -(\mu^2+2\mu t-2t^2)/2,\quad&\text{if}\
0<t+h\le\mu,\\ (6\mu t-3\mu^2-2t^2)/2,\quad&\text{if}\
\mu<t+h\le2\mu,\\ (t-3\mu)^2/2,\quad&\text{if}\ 2\mu<t+h\le3\mu,\\
0,\quad&\text{otherwise}. \end{cases}$$ We  achieve equality in
(\ref{u}) since
$$\|y\|_{L_2}^2=\frac{2}{3}\mu^3,\ \|u \|_{L_2}^2\!=\frac{2}{3\mu^2}
\|z\|_{L_2}^2\!=\frac{2}{3\mu^2}\cdot\mu^5=\frac{2}{3}\mu^3.$$
\end{remark}

\subsection{A Tight Frequency-Domain Stability Criterion} We assume

 {\bf A1} Given the nominal value of the delay $h>0$, the nominal system
\begin{equation}
\begin{array}{lll}
\label{xn}
   \dot x(t)=A_0x(t)+A_1x(t-h),
\end{array}
\end{equation}
is asymptotically stable.

The auxiliary system (\ref{input}) can be written as $y=Gu$ with
the transfer matrix
\begin{equation}
\label{G}
\begin{array}{ll}
G(s)=  \sqrt {{\cal F}(p)} sI (sI-A_0- A_1 e^{-h s})^{-1} \mu A_1
.
\end{array}
\end{equation}
 By the small gain theorem (see
e.g. Gu {\it et al.} (2003)  the system (\ref{x}) is input-output
stable (and thus asymptotically stable, since the nominal system
is time-invariant)
 if
$\|G\|_{\infty}<1$. A stronger result may be obtained by scaling
$G$:
\begin{theorem}
Consider (\ref{x}) with  delay given by (\ref{tau}). Under A1 the
system is asymptotically stable if  there exists non-singular
matrix $X$ such that
\begin{equation}
\label{sg} \|G_X\|_{\infty}<1, \quad G_X(s)= X G(s) X^{-1},
\end{equation}
 where $G$ is given by (\ref{G}) with
${\cal F}(p)$ of (\ref{lem2}) and where $p\in [0, \infty)$
corresponds to case C, while ${\cal F}(\infty)=7/4$ corresponds to
case B.
\end{theorem}

\begin{remark} From Theorem 1 it follows that under A1 (\ref{x}) is
asymptotically stable if
$$\mu< \frac{1}{{\sqrt{{\cal F}(p)}}}\cdot k,
\ k=\frac{1}{\|sI \!(sI\!-\!A_0-\! A_1 \!e^{-h s})^{-1}\!
A_1\|_{\infty}}.$$ By Fridman \& Shaked (2005) ${{\cal F}(p)}=2,
p>0$ and thus (\ref{x}) (with $\dot \tau(t)\leq 1+p, p>0$ or with
$\tau(t)$ of case B) is asymptotically stable for $\tau(t)\in
[h-\mu, h+\mu]$, where $\mu <0.7071k$. By the new bounds of Lemma
\ref{l1} and Lemma \ref{l2} we obtain a wider stability intervals:
\begin{equation}
\label{F}\begin{array}{llllll}
p=0.1,& \dot \tau(t) \leq 1.1, & {{\cal F}(p)}=1.0909, & \mu < 0.9574k,\\
p=0.5, & \dot \tau(t) \leq 1.5, & {{\cal F}(p)}=1.3333, & \mu < 0.8660k,\\
p=1, & \dot \tau(t) \leq 2, & {{\cal F}(p)}=1.5, & \mu < 0.8165 k,\\
p=\infty, &\mbox {case B}, & {{\cal F}(p)}=1.75, & \mu < 0.7559k.
\end{array}\end{equation}
\end{remark}
\subsection { On Improved Stability Criteria in the Time-Domain} By applying the
time-domain results of Fridman \& Shaked (2005)
 via descriptor model transformation with the corresponding
simple Lyapunov-Krasovskii functional
 we
obtain:

\begin{theorem}
System (\ref{x}) is  asymptotically stable  for all delays of
(\ref{tau}), if there exist $n\!\times \!n$ matrices
$0\!<\!P_1,\;P_2,\;P_3, \ S>0$, $Y_{1},\ Y_{2}, \ T,$ $ R
 $, $R_{a}$  such that the following Linear Matrix Inequality
 (LMI)
 \begin{equation}
    \label{LMI2}
{\scriptsize 
\!\left[\begin{array}{ccccccccccc} \Gamma_n &
\begin{array}{ccccccccccc} | \\
| \\
|
\end{array}&
\begin{array}{cc}\mu P_2^TA_1\\
\mu P_3^TA_1\\0
\end{array}&
  \begin{array}{cc}0\\ {\cal F}(p) R_{a}\\ 0 \end{array} \\
   - &  - & - &  -  \\
 {*} & | & -\mu R_{a}&0 \\
 {*} & | & {*}& -{\cal F}(p)R_{a}

    \end{array}\!\!\right]}<0,
    \end{equation}
where
    \begin{equation}\begin{array}{lllll}
\label{Gamma1}
 {
\Gamma_{n}\!=\scriptsize\!\left[\begin{array}{ccccccccccc}
   \Psi_n & P^T\left[\begin{array}{cc}
0 \\
    A_1
 \end{array}\right]-Y^T+\left[\begin{array}{cc} T\\0 \end{array} \right]& hY^T\\
   {*} & -S-T-T^T & -hT\\
   {*} & {*} & -hR
  \end{array}\!\!\right]},\\
 \scriptsize \Psi_{n} = P^T\left[\begin{array}{cc}
0 & I\\
    A_0& -I
 \end{array}\right]+
\left[\begin{array}{cc}
0 & A_0^T\\
    I& -I
 \end{array}\right]P+\left[\begin{array}{cc}S &0\\ 0 &
h R \end{array}\right]\\ \scriptsize+ \left[\begin{array}{cc} Y
\\0
\end{array}\right] +
\left[\begin{array}{cc} Y \\0
\end{array}\right]^T, \ P =\left[ \begin{array}{rr}P_{1} & 0\\ P_2 &P_3
\end{array} \right], \ Y=[Y_{1} \ Y_{2}].
   \end{array} \end{equation}
is feasible. Here ${\cal F}(p)$ is given by (\ref{lem2}) in case C
and ${\cal F}(\infty)=7/4$ in case B.
\end{theorem}

LMI (\ref{LMI2}) is convex in ${\cal F}(p)$ and thus the smaller
values of ${\cal F}(p)$ lead to a less restrictive conditions.
The time-domain criteria give sufficient conditions for the
frequency domain Theorem 1.

{\bf Example} (Kharitonov \& Niculescu, 2003):  Consider the
system
\begin{equation}
\begin{array}{ll}
\label{ex1} \dot x(t)=\left[\begin{array}{cc}0 &1\\ -1
  &-2 \end{array}\right]x(t
  )+\left[\begin{array}{cc}0 &0\\ -1
  & 1 \end{array}\right]x(t-\tau(t)), \\
   \tau(t)=1 +\eta (t), \ |\eta (t)|\leq \mu, \ \dot \tau(t)\leq d.
\end{array}
\end{equation}
 In (Fridman, 2004) the maximum value of $\mu$,
for which the system is asymptotically stable, was found to be
$\mu=0.271$ for all $d\geq 1$. The latter result was less
restrictive than the one by (Kharitonov \& Niculescu, 2003). By
 the time domain criterion of (Fridman \& Shaked, 2005)  for $d=1$
 the corresponding value of $\mu$ is greater ($\mu=0.384$), while for
$d>1$ the result is the same ($\mu =0.271$). Theorem 2 of the
present paper leads to a wider stability interval   for $d>1$ (see
Table 1). Note that the results by Kao \& Rantzer (2005) do not
improve the existing results by descriptor approach and do not
treat the case of $h-\mu >0$.

\begin{table}
\begin{center}
\caption{Maximum value of $\mu$  }
\begin{tabular}{|c|c|c|c|c|c|c|c|c|c|c|} \hline
$d=1$ & $d=1.1$ & $d= 1.5$ & $d=2$ & fast delay  \\
\hline 0.384 & 0.367 &0.331 & 0.313 & 0.289  \\
\hline
\end{tabular}
\end{center}
\end{table}

\section{Conclusions}
Linear systems  with 
bounded time-varying delays are analyzed under the assumption that
the nominal delay values are not equal to zero. Two cases of delay
are considered: case B (without any constraints on the delay
derivative) and case C (where  the delay derivative is not greater
than $d\geq 1$).
 An
input-output approach to stability of such systems  is known to be
based on the bound of the $L_2$-norm of a certain integral
operator. In the present paper for the first time a tight
$d-dependent$ bound is derived.  The existing bound in case B is
also improved. In the past, case C was treated as case B, which
was restrictive.
  The new bounds lead to improved 
stability criteria and gives tools for further improvements.

\section{Appendix}

\begin{figure}
\setlength{\unitlength}{0.5cm}
\begin{picture}(12,9)(0,-2)
\thinlines \put(0.5,1.5){\vector(1,0){10.5}}
\put(0.5,1.5){\vector(0,1){5}} \dashline
{0.2}(0.5,3.5)(11,3.5)\dashline{0.2}(0.5,4.9)(11,4.9)
\put(4.07,2.77){\line(1,1){1.42}}\put(3.84,2.74){\line(1,1){1.46}}\put(3.61,2.71){\line(1,1){1.51}}
\put(3.38,2.68){\line(1,1){1.57}}\put(3.15,2.65){\line(1,1){1.57}}\put(2.92,2.62){\line(1,1){1.54}}
\put(2.69,2.59){\line(1,1){1.54}}\put(2.46,2.56){\line(1,1){1.48}}\put(2.23,2.53){\line(1,1){1.3}}
\put(2,2.5){\line(1,1){1}}
\put(5.75,4.05){\line(1,1){2.14}}\put(5.88,3.98){\line(1,1){2.19}}\put(5.98,3.88){\line(1,1){2.29}}
\put(6.08,3.78){\line(1,1){2.39}}\put(6.21,3.71){\line(1,1){2.46}}\put(6.39,3.69){\line(1,1){2.48}}
\put(6.59,3.69){\line(1,1){2.48}}\put(6.85,3.75){\line(1,1){2.42}}\put(7.35,4.05){\line(1,1){2.12}}
\put(5,3.5){\line(-1,0){2}}\put(5.2,3.7){\line(-1,0){1.9}}\put(5.4,3.9){\line(-1,0){1.75}}
\put(5.6,4.1){\line(-1,0){1.55}}\put(5.6,4.1){\line(1,0){1.7}}\put(5.9,3.9){\line(1,0){1.2}}
\put(5.8,4.3){\line(1,0){1.85}}\put(6,4.5){\line(1,0){1.85}}\put(6.2,4.7){\line(1,0){1.85}}
\put(6.4,4.9){\line(1,0){1.85}} \thicklines
\put(1.8,2.8){\line(1,1){3.7}}\put(4.3,2.8){\line(1,1){3.6}}\put(6.8,2.8){\line(1,1){3.7}}
\qbezier(1.7,2.3)(3.5,4.2)(5,4.3)\qbezier(5,4.3)(5.7,4.2)(6,3.8)\qbezier(6,3.8)(6.5,3.4)(7.5,4.2)
\qbezier(7.5,4.2)(8.5,5.1)(9.6,6.2)
\put(11,1.7){$s_1$}\put(0.7,6.5){$s_2$}\put(-0.1,3.4){$t_1$}\put(-0.1,4.8){$t_2$}
\put(0,0.5){$t_1=t-h-\eta(t),\quad
t_2=t-h$}\put(0,-0.5){$\Lambda_1: s_2=s_1-h+\mu,\quad \Lambda_2:
s_2=s_1-h,$}\put(0,-1.5){$\Lambda_3: s_2=s_1-h-\mu,\quad \Gamma:
s_2=s_1-h-\eta(s_1)$}\put(5,6.8){$\Lambda_1$}\put(7.9,6.6){$\Lambda_2$}\put(10.5,6.8){$\Lambda_3$}
\put(1.1,1.8){$\Gamma$}
\end{picture}
\caption{Domains $D$ and $D(t)$}\label{figfs4}
\end{figure}
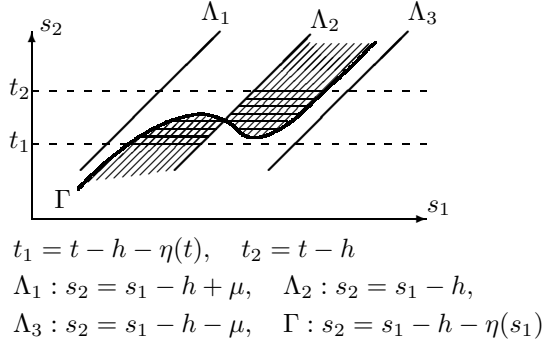

{\bf Proof of Lemma \ref{l1}}. Denote by $\varphi:\R^2\to\{0,1\}$
the characteristic function of the domain $D$ in the positive
quadrant, bounded by the line $s=t-h$ and by the graph of the
function $s=t-h-\eta(t)$, i.e.,
$$\varphi(t,s)\!=\!\begin{cases}1,\ &\text{if}\ (t-h-s)(t-h-\eta(t)-s)\le
0,\\ 0,\ &\text{if}\ (t-h-s)(t-h-\eta(t)-s)>0\end{cases}$$ (shaded
region in Figure \ref{figfs4}). Then $z(t)$ given by (\ref{Del})
satisfies the following:
$$\begin{array}{lllllll}\|z\|_{L_2}^2=\int_0^{\infty}
\left(\int_{-\infty}^{\infty}\varphi(t,s_1)y(s_1)ds_1\right)^T\\
\times
\left(\int_{-\infty}^{\infty}\varphi(t,s_2)y(s_2)ds_2\right)dt\\
=\int_{-\infty}^{\infty}\int_{-\infty}^{\infty}\left(\int_0^{\infty}
\varphi(t,s_1)\varphi(t,s_2)dt\right)\\
 \times
y^T(s_1)y(s_2)ds_1ds_2\\
=\int_{-\infty}^{\infty}\int_{-\infty}^{\infty}k(s_1,s_2)y^T(s_1)y(s_2)ds_1ds_2,\\
k(s_1,s_2)= \int_0^{\infty}\varphi(t,s_1)\varphi(t,s_2)dt,\
s_1,s_2\in\R\ .\end{array}$$
Hence $\|z\|_{L_2}^2\le\|{\cal
K}\|_{L_2}\cdot\|y\|_{L_2}^2$, where $\|{\cal K}\|_{L_2}$ is the
$L_2$-norm of the 
operator ${\cal K}:L_2(\R)\to L_2(\R)$
$${\cal K}(f)(t)=\int_{-\infty}^\infty k(t,s)f(s)ds,\quad f\in
L_2(\R)\ .$$
By the Riesz-Thorin interpolation theorem (see, for example,
Okikiolu (1971), Theorem 5.1.3), 
$\|{\cal K}\|_{L_2}\le\sqrt{\|{\cal
K}\|_{L_1}\cdot\|{\cal K}\|_{L_\infty}}$. Since $k(s_1,s_2)\ge 0$
and $k(s_1,s_2)=k(s_2,s_1)$, by the well-known formulas for the
$L_1$ and $L_\infty$-norms, we have $\|{\cal K}\|_{L_1}=\|{\cal
K}\|_{L_\infty}\le\sup_{s\in[0,\infty)}K(s)$, where
$K(s)=\int_0^\infty k(s_1,s)ds_1$, and hence we decide that
\begin{equation}\|{\cal K}\|_{L_2}\le\sup_{s\in[0,\infty)}K(s)\ .\label{e1}\end{equation}
We shall show that $K(t)\le7/4\ \mu^2$ for all $t\in[0,\infty)$.

Without loss of generality assume that $\eta(t)>0$. Geometrically,
$K(t)$ is the area of the part $D(t)$ of the domain $D$ cut out by
the strip $t-h\ge s_2\ge t-h-\eta(t)$ (double shaded region in
Figure \ref{figfs4}).

\begin{figure}
\setlength{\unitlength}{0.5cm}
\begin{picture}(13,7.5)(0,-6.5)
\put(6,1.2){$\Lambda_1$}\put(8.3,0.7){$\Lambda_2$}\put(11.7,1.2){$\Lambda_3$}
\put(0,-5){$t_3=t-\mu-\eta(t),\ t_4=t-\eta(t),\
t_5=t,$}\put(1,-6){$t_6=t+\mu,\ t_7=(t_4+t_5)/2$} \thinlines
\put(2,-3){\vector(1,0){11}}
\put(2,-3){\vector(0,1){4}}\put(7,-2){\line(0,1){2}}\put(3.2,-2){\line(1,1){2}}\put(3.4,-2){\line(1,1){2}}
\put(3.6,-2){\line(1,1){2}}\put(3.8,-2){\line(1,1){2}}\put(4,-2){\line(1,1){2}}\put(4.2,-2){\line(1,1){2}}
\put(4.4,-2){\line(1,1){2}}\put(4.6,-2){\line(1,1){2}}\put(4.8,-2){\line(1,1){2}}\put(5,-2){\line(1,1){2}}
\put(5.2,-2){\line(1,1){1.8}}\put(5.4,-2){\line(1,1){1.6}}\put(5.6,-2){\line(1,1){1.4}}
\put(5.8,-2){\line(1,1){1.2}}\put(6,-2){\line(1,1){2}}\put(7,-1.2){\line(1,1){1.2}}\put(7,-1.4){\line(1,1){1.4}}
\put(7,-1.6){\line(1,1){1.6}}\put(7,-1.8){\line(1,1){1.8}}\put(7,-2){\line(1,1){2}}\put(7.2,-2){\line(1,1){2}}
\put(7.4,-2){\line(1,1){2}}\put(7.6,-2){\line(1,1){2}}\put(7.8,-2){\line(1,1){2}}\put(8,-2){\line(1,1){2}}
\put(8.2,-2){\line(1,1){2}}\put(8.4,-2){\line(1,1){2}}\put(8.6,-2){\line(1,1){2}}\put(8.8,-2){\line(1,1){2}}
\dashline{0.2}(2,-2)(3,-2)\dashline{0.2}(2,0)(5,0)
\dashline{0.2}(2.3,-2.7)(3,-2)\dashline{0.2}(5,0)(6,1)\dashline{0.2}(5.3,-2.7)(8.5,0.5)\dashline{0.2}(8.3,-2.7)(9,-2)
\dashline{0.2}(11,0)(12,1)
\dottedline{0.1}(3,-3)(3,-2)\dottedline{0.1}(6,-3)(6,0)\dottedline{0.1}(7,-3)(7,-2)\dottedline{0.1}(11,-3)(11,0)
\dottedline{0.1}(8,-3)(8,0) \thicklines
\put(3,-2){\line(1,1){2}}\put(9,-2){\line(1,1){2}}\put(3,-2){\line(1,0){6}}\put(5,0){\line(1,0){6}}
\put(2.9,-3.6){$t_3$}\put(5.9,-3.6){$t_4$}\put(6.9,-3.6){$t_7$}\put(7.9,-3.6){$t_5$}\put(10.9,-3.6){$t_6$}
\put(1.4,-2.2){$t_1$}\put(1.4,-0.2){$t_2$}\put(12.8,-2.8){$s_1$}\put(2.3,1){$s_2$}
\end{picture}
\caption{Parallelogram $\Pi(t)$ and domain $D(t)$}\label{figfs1}
\end{figure}
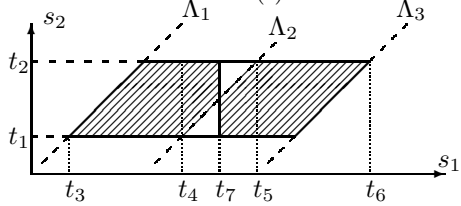

Thus, $D(t)$ lies inside the parallelogram
$$\begin{array}{lll}\Pi(t):=\{(s_1,s_2)\in\R^2 \!: t-h\le s_2\le t-h-\eta(t),\\
s_1-h-\mu\le s_2\le s_1-h+\mu\}\end{array}$$ (see Figure
\ref{figfs1}). Divide $\Pi(t)$ by the vertical lines $s_1=t_4$ and
$s_1=t_5$ into two trapezes of total area $\mu^2-(\mu-\eta(t))^2$,
and a square (see Figure \ref{figfs1}). The intersection of $D(t)$
with a vertical line $s_1=\sigma$, where $t_4\le\sigma\le t_5$, is
contained either in the segment $[t-h-\eta(t),\tau-h]$, or in the
segment $[\tau-h,t-h]$, and hence the area of
$D(t)\cap\{t-\eta(t)\le s_1\le t\}$ does not exceed
$$\int_{t-\eta(t)}^t\max\{(\tau-(t-\eta(t)),t-\tau\}d\tau=\frac{3}{4}\eta(t)^2\
.$$ So we derive the required bound (\ref{en1}) from the evident
inequality
$$\mu^2-(\mu-\eta(t))^2+\frac{3}{4}\eta(t)^2\le\frac{7}{4}\mu^2\
,$$which geometrically means that the maximal area domain $D(t)$
looks as the shaded region in Figure \ref{figfs1}.

\begin{figure}
\setlength{\unitlength}{0.5cm}
\begin{picture}(16,10)(0,-5)
\put(3.5,4.5){$p\ge1$}\put(11.5,4.5){$p<1$} \thicklines
\put(0.5,1.5){\vector(1,0){6}}
\put(0.5,1.5){\vector(0,1){3}}\put(8.5,1.5){\vector(1,0){6}}
\put(8.5,1.5){\vector(0,1){3}}\put(0.5,1.5){\line(1,1){2}}\put(4.5,1.5){\line(1,1){2}}
\put(8.5,1.5){\line(1,1){2}}\put(12.5,1.5){\line(1,1){2}}\put(2.5,3.5){\line(1,0){4}}
\put(10.5,3.5){\line(1,0){4}}\put(4,1.5){\line(-1,2){1}}\put(12.5,1.5){\line(-3,1){3}}
\put(2.5,1.5){\line(1,1){2}}\put(10.5,1.5){\line(1,1){2}}
\thinlines\dashline{0.2}(0.5,3.5)(2.5,3.5)\dashline{0.2}(8.5,3.5)(10.5,3.5)
\dashline{0.2}(2.5,1.5)(2.5,3.5)\dashline{0.2}(10.5,1.5)(10.5,3.5)\dashline{0.2}(4.5,1.5)(4.5,3.5)
\dashline{0.2}(12.5,1.5)(12.5,3.5)\put(0.7,1.5){\line(1,1){2}}\put(0.9,1.5){\line(1,1){2}}
\dashline{0.2}(3.5,1.5)(3.5,2.5)
\put(1.1,1.5){\line(1,1){1.933}}\put(1.3,1.5){\line(1,1){1.8}}\put(1.5,1.5){\line(1,1){1.667}}
\put(1.7,1.5){\line(1,1){1.533}}
\put(1.9,1.5){\line(1,1){1.4}}\put(2.1,1.5){\line(1,1){1.267}}\put(2.3,1.5){\line(1,1){1.133}}
\put(4.7,3.5){\line(-1,-1){1.133}}\put(4.9,3.5){\line(-1,-1){1.267}}\put(5.1,3.5){\line(-1,-1){1.4}}
\put(5.3,3.5){\line(-1,-1){1.533}}\put(5.5,3.5){\line(-1,-1){1.667}}\put(5.7,3.5){\line(-1,-1){1.8}}
\put(5.9,3.5){\line(-1,-1){1.933}}\put(6.1,3.5){\line(-1,-1){2}}\put(6.3,3.5){\line(-1,-1){2}}
\put(8.7,1.5){\line(1,1){0.95}}\put(8.9,1.5){\line(1,1){0.9}}\put(9.1,1.5){\line(1,1){0.85}}
\put(9.3,1.5){\line(1,1){0.8}}\put(9.5,1.5){\line(1,1){0.75}}\put(9.7,1.5){\line(1,1){0.7}}
\put(9.9,1.5){\line(1,1){0.65}}\put(10.1,1.5){\line(1,1){0.6}}\put(10.3,1.5){\line(1,1){0.55}}
\put(10.5,1.5){\line(1,1){0.5}}\put(12.7,3.5){\line(-1,-1){1.55}}\put(12.9,3.5){\line(-1,-1){1.6}}
\put(13.1,3.5){\line(-1,-1){1.65}}\put(13.3,3.5){\line(-1,-1){1.7}}\put(13.5,3.5){\line(-1,-1){1.75}}
\put(13.7,3.5){\line(-1,-1){1.8}}\put(13.9,3.5){\line(-1,-1){1.85}}\put(14.1,3.5){\line(-1,-1){1.9}}
\put(14.3,3.5){\line(-1,-1){1.95}}\put(2.2,0.8){$\mu$}\put(3.2,0.8){$\frac{3}{2}\mu$}\put(4.4,0.8){$2\mu$}
\put(0,3.4){$\mu$}
\put(10.4,0.8){$\mu$}\put(12.3,0.8){$2\mu$}\put(8,3.4){$\mu$}
\put(3.5,-0.1){{\rm (a)}}\put(11.5,-0.1){{\rm (b)}}
\thicklines\put(4,-4){\line(1,1){2}}\put(10,-4){\line(1,1){2}}\put(4,-4){\line(1,0){6}}
\put(6,-2){\line(1,0){6}}\thinlines\put(6.7,-4.3){\line(1,1){3.3}}\put(9.1,-4.2){\line(-1,2){1.6}}
\put(9,-4.3){\line(0,1){3.3}}\put(8.1,-4.2){\line(-1,2){1.6}}
\put(4.2,-4){\line(1,1){2}}\put(4.4,-4){\line(1,1){2}}\put(4.6,-4){\line(1,1){2}}
\put(4.8,-4){\line(1,1){2}}\put(5,-4){\line(1,1){2}}\put(5.2,-4){\line(1,1){1.867}}
\put(5.4,-4){\line(1,1){1.733}}\put(5.6,-4){\line(1,1){1.6}}\put(5.8,-4){\line(1,1){1.467}}
\put(6,-4){\line(1,1){1.333}}\put(6.2,-4){\line(1,1){1.2}}\put(6.4,-4){\line(1,1){1.067}}
\put(6.6,-4){\line(1,1){0.933}}\put(6.8,-4){\line(1,1){0.8}}\put(9.2,-2){\line(-1,-1){1.467}}
\put(9.4,-2){\line(-1,-1){1.6}}\put(9.6,-2){\line(-1,-1){1.733}}\put(9.8,-2){\line(-1,-1){1.867}}
\put(10,-2){\line(-1,-1){2}}\put(10.2,-2){\line(-1,-1){2}}\put(10.4,-2){\line(-1,-1){2}}
\put(10.6,-2){\line(-1,-1){2}}\put(10.8,-2){\line(-1,-1){2}}\put(11,-2){\line(-1,-1){2}}
\put(11.2,-2){\line(-1,-1){2}}\put(11.4,-2){\line(-1,-1){2}}\put(11.6,-2){\line(-1,-1){2}}
\put(11.8,-2){\line(-1,-1){2}}\put(6.5,-0.9){$l_3$}\put(7.5,-0.9){$l_2$}\put(9,-0.9){$l_1$}
\put(10.1,-0.9){$\Lambda_2$}\put(7.9,-5){{\rm (c)}}
\end{picture}
\caption{The case $N=1$}\label{figfs2}
\end{figure}
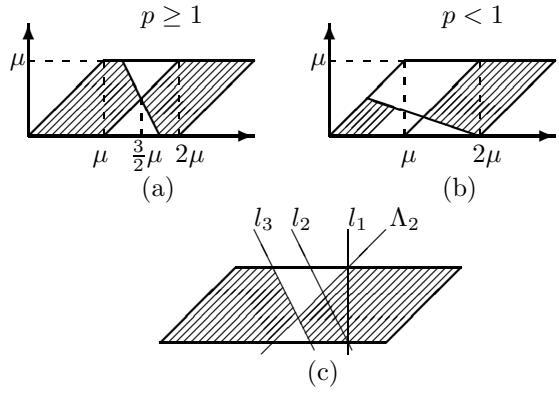

{\bf Proof of Lemma \ref{l2}}. We interpret the function ${\cal
F}(p)$ geometrically as follows:

- for $p\ge 1$, ${\cal F}(p)\mu^2$ is the area of the domain
$$\begin{array}{lll} D_p:=\{(t,s)\in\R^2\ :\ 0\le s\le\mu,\ t-2\mu\le s\le
t,\\
(t-\mu-s)(2s+2pt-(3p+1)\mu))\ge 0\}\end{array}$$ (shaded region in
Figure \ref{figfs2}(a)),

- for $0\le p<1$, ${\cal F} (p)\mu^2$ is the area of the domain
$$\begin{array}{lll} D_p=\{(t,s)\in\R^2\ :\ 0\le s\le \mu,\ t-2\mu\le s\le
t,\\
(t-\mu-s)(s+pt-2p\mu)\le 0\}\end{array}$$ (shaded region in Figure
\ref{figfs2}(b)).

As in the proof of Lemma \ref{l1}, we estimate from above the
value of $K(s_2)$. i.e., the area of the domain $D(t)$.

Without loss of generality we assume that $\eta$ is a smooth
function, whose zero locus is locally finite and consists of only
simple roots. Fix $t>0$. Assume that $\eta(t)>0$, which, in
particular, means that the line $s_1=t$ crosses the domain $D(t)$
along the segment $t_1\le s_2\le t_2$ (cf. Figure \ref{figfs1}).
Put $N=\#(\eta^{-1}(0)\cap(t-\mu-\eta(t),t+\mu))$ (i.e., the
number of zeroes in the interval $(t_3,t_6)$ in Figure
\ref{figfs1}).

Consider few possibilities.

{\it Step 1}. Suppose that $N=0$. Then $D(t)$ is contained in the
part of the parallelogram $\Pi(t)$, right to the line $s_2=s_1-h$
(see Figure \ref{figfs1}). Then its area does not exceed
$\eta(t)\mu\le\mu^2\le{\cal F}(p)\mu^2$.

{\it Step 2}. Suppose that $N=1$. The zero
$\tau=\eta^{-1}(0)\cap(t-\mu-\eta(t),t+\mu)$ corresponds to an
intersection point of the graphs of the functions $s_2=s_1-h$ and
$s_2=s_1-h-\eta(s_1)$. This intersection point lies either right
to the line $l_1:=\{s_1=t\}$ or below the line
\mbox{$l_2:=\{s_2=-ps_1+(p-1)t-h-\eta(t)\}$} (see Figure
\ref{figfs2}(c)). In the former case, the domain $D(t)$ remains
below the line $s_2=s_1-h$, that is we have the upper bound from
Step 1. In the latter case, the part of domain $D(t)$ in the
half-plane $s_1\ge\tau$ should lie below the line $s_2=s_1-h$ and
above the line $l_3=\{s_2=-ps_1+(p-1)\tau-h\}$, and the part of
$D(t)$ in the half-plane $s_1\le\tau$ should lie below the line
$l_3$ and above the line $s_2=s_1-h$ (shaded region in Figure
\ref{figfs2}(c)). It is an elementary geometry exercise to show
that the area of the shaded region in Figure \ref{figfs2}(c) does
not exceed ${\cal F}(p)\mu^2$.

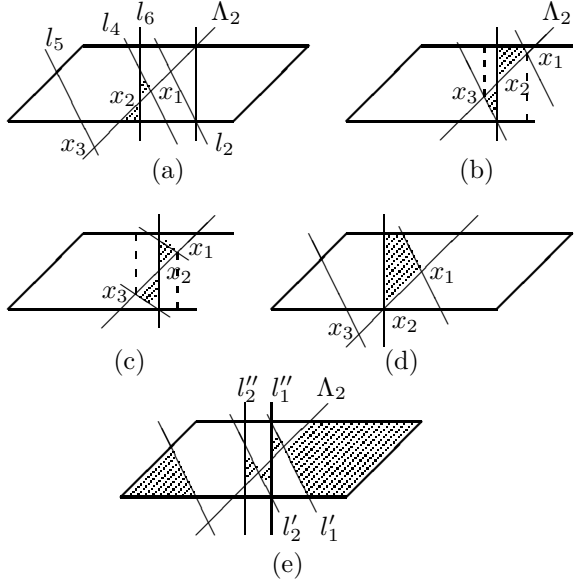
\begin{figure}
\setlength{\unitlength}{0.5cm}
\begin{picture}(16,15)(0,0)
\thicklines\put(3,2){\line(1,0){6}}\put(5,4){\line(1,0){6}}
\put(3,2){\line(1,1){2}}\put(9,2){\line(1,1){2}}\put(0,7){\line(1,0){5}}
\put(0,7){\line(1,1){2}}\put(2,9){\line(1,0){4}}\put(7,7){\line(1,0){6}}
\put(7,7){\line(1,1){2}}\put(9,9){\line(1,0){6}}\put(13,7){\line(1,1){2}}
\put(0,12){\line(1,0){6}}\put(0,12){\line(1,1){2}}\put(6,12){\line(1,1){2}}
\put(2,14){\line(1,0){6}}\put(9,12){\line(1,0){5}}\put(9,12){\line(1,1){2}}
\put(11,14){\line(1,0){4}}
\thinlines\put(3.5,11.5){\line(0,1){3}}\put(2,11){\line(1,1){3.5}}\put(4.5,11.4){\line(-1,2){1.4}}
\put(2.4,11.1){\line(-1,2){1.4}}
\put(5.3,11.4){\line(-1,2){1.4}}\put(5,11.5){\line(0,1){3}}
\put(11.5,11.5){\line(1,1){3}}\put(13.25,11.5){\line(-1,2){1.3}}
\put(13,11.5){\line(0,1){3}}\put(14.5,12.4){\line(-1,2){1}}\put(2.5,6.5){\line(1,1){3}}
\put(4.3,6.8){\line(-3,2){1.3}}\put(4,6.5){\line(0,1){3}}\put(4.8,8.3){\line(-3,2){1.3}}
\put(9,6){\line(1,1){3.5}}\put(9.5,6){\line(-1,2){1.5}}\put(10,6){\line(0,1){3.5}}
\put(11.7,6.6){\line(-1,2){1.3}}\put(5,1){\line(1,1){3.5}}\put(5.5,1){\line(-1,2){1.5}}
\put(6.3,1){\line(0,1){3.5}}\put(7.2,1.6){\line(-1,2){1.3}}\put(7,1){\line(0,1){3.5}}
\put(8.2,1.6){\line(-1,2){1.3}}
\dashline{0.2}(12.667,12.667)(12.667,14)\dashline{0.2}(13.8,13.8)(13.8,12)
\dashline{0.2}(3.4,7.4)(3.4,9)\dashline{0.2}(4.5,8.5)(4.5,7)
\dottedline{0.1}(3.5,12.2)(3.3,12)\dottedline{0.1}(3.5,12.4)(3.1,12)\dottedline{0.1}(3.5,12.7)(3.733,12.933)
\dottedline{0.1}(3.5,12.9)(3.633,13.033)
\dottedline{0.1}(13,12.8)(12.733,12.533)\dottedline{0.1}(13,12.6)(12.8,12.4)\dottedline{0.1}(13,12.4)(12.867,12.267)
\dottedline{0.1}(13,13.2)(13.72,13.92)\dottedline{0.1}(13,13.4)(13.6,14)
\dottedline{0.1}(13,13.6)(13.4,14)\dottedline{0.1}(13,13.8)(13.2,14)
\dottedline{0.1}(4,7.4)(3.75,7.2)\dottedline{0.1}(4,7.6)(3.64,7.24)\dottedline{0.1}(4,7.8)(3.52,7.32)
\dottedline{0.1}(4,8.2)(4.4,8.6)\dottedline{0.1}(4,8.4)(4.3,8.7)\dottedline{0.1}(4,8.6)(4.2,8.8)
\dottedline{0.1}(10,7.2)(10.933,8.133)\dottedline{0.1}(10,7.4)(10.867,8.267)
\dottedline{0.1}(10,7.6)(10.8,8.4)\dottedline{0.1}(10,7.8)(10.733,8.533)
\dottedline{0.1}(10,8)(10.667,8.667)\dottedline{0.1}(10,8.2)(10.6,8.8)
\dottedline{0.1}(10,8.4)(10.533,8.933)\dottedline{0.1}(10,8.6)(10.4,9)
\dottedline{0.1}(10,8.8)(10.2,9)\dottedline{0.1}(3.2,2)(4.4,3.2)\dottedline{0.1}(3.4,2)(4.4067,3.033)
\dottedline{0.1}(3.6,2)(4.533,2.933)\dottedline{0.1}(3.8,2)(4.6,2.8)\dottedline{0.1}(4,2)(4.667,2.667)
\dottedline{0.1}(4.2,2)(4.733,2.533)\dottedline{0.1}(4.4,2)(4.8,2.4)
\dottedline{0.1}(4.6,2)(4.867,2.267)\dottedline{0.1}(4.8,2)(4.933,2.133)
\dottedline{0.1}(6.3,2.5)(6.6,2.8)\dottedline{0.1}(6.3,2.7)(6.5,2.9)\dottedline{0.1}(6.3,2.9)(6.4,3)
\dottedline{0.1}(7,2.8)(6.733,2.533)\dottedline{0.1}(7,2.6)(6.867,2.467)
\dottedline{0.1}(8.2,4)(7.4,3.2)\dottedline{0.1}(8.4,4)(7.467,3.067)
\dottedline{0.1}(8.6,4)(7.533,2.933)\dottedline{0.1}(8.8,4)(7.6,2.8)\dottedline{0.1}(9,4)(7.667,2.677)
\dottedline{0.1}(9.2,4)(7.733,2.533)\dottedline{0.1}(9.4,4)(7.8,2.4)
\dottedline{0.1}(9.6,4)(7.867,2.267)\dottedline{0.1}(9.8,4)(7.933,2.133)
\dottedline{0.1}(10,4)(8,2)\dottedline{0.1}(10.2,4)(8.2,2)\dottedline{0.1}(10.4,4)(8.4,2)
\dottedline{0.1}(10.6,4)(8.6,2)\dottedline{0.1}(10.8,4)(8.8,2)
\dottedline{0.1}(7,3.2)(7.267,3.467)\dottedline{0.1}(7,3.4)(7.133,3.533)
\put(3.8,10.5){{\rm (a)}}\put(12,10.5){{\rm
(b)}}\put(2.8,5.4){{\rm (c)}}\put(10,5.4){{\rm (d)}}\put(7,0){{\rm
(e)}}\put(1.4,11.2){$x_3$}\put(2.7,12.5){$x_2$}\put(4,12.6){$x_1$}\put(5.4,14.7){$\Lambda_2$}
\put(14.2,14.7){$\Lambda_2$}
\put(5.5,11.2){$l_2$}\put(2.5,14.4){$l_4$}\put(1,14){$l_5$}\put(3.4,14.7){$l_6$}
\put(14.1,13.5){$x_1$}\put(13.2,12.8){$x_2$}\put(11.8,12.6){$x_3$}
\put(4.8,8.45){$x_1$}\put(4.2,7.8){$x_2$}\put(2.5,7.3){$x_3$}
\put(11.2,7.8){$x_1$}\put(10.2,6.5){$x_2$}\put(8.5,6.3){$x_3$}\put(8.2,4.7){$\Lambda_2$}
\put(8.3,1){$l'_1$}\put(7.3,1){$l'_2$}\put(7,4.7){$l''_1$}\put(6.1,4.7){$l''_2$}
\end{picture}
\caption{The case $N>1$}\label{figfs3}
\end{figure}

{\it Step 3}. We intend to show that the case $N>1$ reduces to the
above considered cases $N=0$ or $1$. For, we need the following
auxiliary geometric statement. Consider the parallelogram
$\Pi(t)$, some points $x_1,x_2,x_3$ with decreasing coordinates,
lying on the line $s_2=s_1-h$ below the line
$l_2=\{s_2=-ps_1+(p-1)t-h-\eta(t)\}$ (see Figure \ref{figfs3}(a)).
Draw the lines $l_4,l_5$ with slope $-p$ through the points
$x_1,x_3$, respectively, and the vertical line $l_6$ through
$x_2$. Denote by $F(x_2)$ the area of the domain $\del(x_2)$,
lying inside $\Pi(t)$, between the lines $l_4,l_5$, and in the two
sectors, bounded by the lines $l_6$ and $s_2=s_1-h$ as marked by
dots in Figure \ref{figfs3}(a). Consider $F(x_2)$ as a function of
the point $x_2$ running along the segment $[x_1,x_3]$. We claim
that $F$ attains its maximum either at $x_2=x_1$ or at $x_2=x_3$.

Indeed, when $x_2$ runs along a subsegment of $[x_1,x_3]$ such
that the combinatorial type of $\del(x_2)$ does not change, $F$ is
a quadratic polynomial in the abscissa $\tau$ of $x_2$ with the
positive second derivative. That is, in this subsegment, $F$
attains its maximum at an endpoint. If such an endpoint differs
from $x_1$ and $x_3$, then it corresponds to the situation when
\begin{itemize}\item either $x_2$ belongs to a side of $\Pi(t)$ (see, for example, Figure
\ref{figfs3}(d)); \item or the domain $\del(x_2)$ intersects with
some side of $\Pi(t)$ at a point (see, Figure
\ref{figfs3}(b,c)).\end{itemize} In the former situation,
$\del(x_2)$ entirely lies right to the vertical line through
$x_2$, and then monotonically grows as $x_2$ tends to $x_3$. In
the latter situation, when replacing $x_2$ by $x_1$, the domain
$\del(x_2)$ turns into the domain $\del(x_1)$ by getting rid of
$\del(x_2)\cap\{s_1\le\tau\}$ and adding the fragment $\del_-$
(the trapeze, bounded by the vertical line through $x_2$, dashed
line, the upper side of $\Pi(t)$, and the line $s_2=s_1-h$ in
Figure \ref{figfs3}(b)), and when replacing $x_2$ by $x_3$, the
domain $\del(x_2)$ turns into the domain $\del(x_3)$ by by getting
rid of $\del(x_2)\cap\{s_1\ge\tau\}$ and adding the fragment
$\del_+$ (the trapeze, bounded by the vertical line through $x_2$,
dashed line, the lower side of $\Pi(t)$, and the line $s_2=s_1-h$
in Figure \ref{figfs3}(b)). One can easy see that the area of
$\del(x_2)$ in the above situation is less than the maximum of the
areas of $\del(x_1)$ and $\del(x_3)$. The same conclusion one can
derive in the situation, presented in Figure \ref{figfs3}(c).

{\it Step 4}. Suppose that $N>1$. We can take $N$ to be odd,
adding if necessary one more zero close to $t-\mu-\eta(t)$. Denote
by $x_1,x_2,...,x_n$ the corresponding intersection points of the
graph of $s_2=s_1-h-\eta(s_1)$ with the line $s_2=s_1-h$, numbered
by the decreasing coordinates. Through each point $x_{2i-1}$ we
draw a line $l'_i$ with slope $-p$, $1\le i\le(n+1)/2$, and
through each point $x_{2i}$ we draw a vertical line $l''_i$, $1\le
i\le n/2$. Thus, $D(t)$ is the union of the regions in $\Pi(t)$,
bounded by the introduced lines as follows (marked by dots in
Figure \ref{figfs3}(e)):
\begin{itemize}\item above the line $l'_1$ and below the line
$s_2=s_1-h$, \item below the line $l'_i$, right to the line
$l''_i$, and above the line $s_2=s_1-h$, $1\le i\le n/2$, \item
left to the line $l''_i$, above the line $l'_{i+1}$ , and below
the line $s_2=s_1-h$, $1\le i\le n/2$, \item below the line
$l'_{(n+1)/2}$ and above the line $s_2=s_1-h$.
\end{itemize}

Using the statement of Step 3, we can move the point $x_2$ either
to the position $x_1$, or $x_3$ and increase the area of
$\del(t)$. On the other hand, each of these limit positions for
$x_2$ means that we, in fact, have reduced two zeroes of $\eta$ in
the interval $(t-\mu-\eta(t),t+\mu)$. Thus, we inductively come to
the case $N=1$ treated in Step 2.

\begin{center}
\begin{Large}
{\bf References}
\end{Large}
\end{center}
\begin{quote}
\hspace*{-1cm}
 Fridman, E.
 (2004).
 Stability of linear functional
differential equations: A new Lyapunov technique. In {\it
Proceedings of MTNS 2004, Leuven}.
\\
\hspace*{-1cm}
%
 Fridman, E. \& Shaked, U. (2005).  Stability and $L_2-$Gain
Analysis of Systems with Time-Varying Delays: Input-Output
Approach.  In:
 {\em  Proc. of 44th  Conf. on Decision and
 Control}, Sevilla, Spain, 2005.
 \\
\hspace*{-1cm}
 Gu, K.,  Kharitonov, V. \& J. Chen, J. (2003).  {\it Stability of
time-delay systems}.  Birkhauser: Boston.
\\
\hspace*{-1cm}
 Kao, C.-Y. \&  Lincoln, B. (2004). Simple stability criteria for systems
with time-varying delays. {\it Automatica}, 40, 1429-1434.
\\
\hspace*{-1cm} Kao, C.-Y. \& Rantzer, A. (2005) Robust stability
analysis of linear systems with time-varying delays. In: {\em
Proc. of 16-th IFAC World Congress}, Prague, Chesh Republic, July
2005.
\\
\hspace*{-1cm}
Kharitonov, V. \&  Niculescu, S. (2003). On the stability of
linear systems with uncertain delay. {\it IEEE Trans. on Automat.
Contr.}, 48, 127-132.
\\
\hspace*{-1cm}
 Kolmanovskii, V. \& Myshkis, A. (1999). {\it
Applied Theory of functional differential equations}, Kluwer.
\\
\hspace*{-1cm}
 Niculescu, S.-I. (2001).
{ Delay effects on stability: A Robust Control Approach}, {\it
 Lecture Notes in
Control and Information Sciences}, {269}, Springer-Verlag, London.
\\
\hspace*{-1cm}
Okikiolu, G. O. (1971). {\it Aspects of the Theory of Bounded
Integral Operators in $L^p$-spaces}, Acad. Press, London.

\end{quote}

\end{document}